\documentclass[11pt]{article}
\usepackage{amsmath,amssymb}
\usepackage{amsthm} 
\usepackage{latexsym}
\usepackage{mathtools}
\newtheorem{thm}{Theorem}[section]

\newtheorem{lem}[thm]{Lemma}
\newtheorem{cor}[thm]{Corollary}

\newcommand{\bpf}{\noindent {\bf Proof}. $\mbox{}$}
\newcommand{\epf}{\hfill \hbox{\raisebox{.5ex}{\fbox}$\phantom{.}$}}

\begin{document}
\title{\bf The Facets of the Subtours Elimination Polytope}
\date{}
\maketitle
\begin{center}
\author{
{\bf Brahim Chaourar} \\
{ Department of Mathematics and Statistics,\\Al Imam Mohammad Ibn Saud Islamic University (IMSIU) \\P.O. Box
90950, Riyadh 11623,  Saudi Arabia }\\{Correspondence address: P. O. Box 287574, Riyadh 11323, Saudi Arabia}}
\end{center}



\begin{abstract}
\noindent Let $G=(V, E)$ be an undirected graph. The subtours elimination polytope $P(G)$ is the set of $x\in \mathbb{R}^E$ such that: $0\leq x(e)\leq 1$ for any edge $e\in E$, $x(\delta (v))=2$ for any vertex $v\in V$, and $x(\delta (U))\geq 2$ for any nonempty and proper subset $U$ of $V$. $P(G)$ is a relaxation of the Traveling Salesman Polytope, i.e., the convex hull of the Hamiton circuits of $G$. Maurras \cite{Maurras 1975} and Gr\"{o}tschel and Padberg \cite{Grotschel and Padberg 1979b} characterize the facets of $P(G)$ when $G$ is a complete graph. In this paper we generalize their result by giving a minimal description of $P(G)$ in the general case and by presenting a short proof of it.

\end{abstract}

\noindent {\bf2010 Mathematics Subject Classification:} Primary 90C57, Secondary 90C27, 52B40. \newline {\bf Key words and phrases:} Traveling Salesman Problem, subtours elimination constraints, facets, locked subgraphs.
\section{Introduction}

Sets and their characteristic vectors will not be distinguished. We refer to Bondy and Murty \cite{Bondy and Murty 2008}, Oxley \cite{Oxley 1992} and Schrijver \cite{Schrijver 2004} about, respectively, graphs, matroids and polyhedra terminolgy and facts.

\noindent Let $G=(V, E)$ be an undirected graph. The subtours elimination polytope $P(G)$ is the set of $x\in \mathbb{R}^E$ such that:
$$  x(e)\geq 0  \> \mathrm{for\>any\> edge}\> e\in E                                                \eqno    (1)$$
$$	x(e) \leq 1 \> \mathrm{for\>any\> edge}\> e\in E			                                    \eqno	 (2)$$
$$	x(\delta (v))=2 \> \mathrm{for\>any\> vertex}\> v\in V			                                \eqno	 (3)$$
$$	x(\delta (U))\geq 2	 \> \mathrm{for\> any\> nonempty\> and\> proper\> subset}\> U\subseteq V 	\eqno	 (4)$$

\noindent Inequalities (1) and (2) are called the trivial constraints, equalities (3) the degree constraints, and inequalities (4) the subtours elimination constraints. $P(G)$ is a relaxation of the Traveling Salesman Polytope (TSP), i.e., the convex hull of the Hamiton circuits of $G$, which permits to get a lower bound for the TSP \cite{Dantzig et al. 1954}. Maurras \cite{Maurras 1975} and Gr\"{o}tschel and Padberg \cite{Grotschel and Padberg 1979b} characterize when a subtour elimination constraint is a facet of $P(K_n)$ with $n\geq 3$ as follows.

\begin{thm}
An inequality (4) defines a facet of $P(K_n)$ ($n\geq 3$) if and only if $2\leq |U|\leq n-2$.
\end{thm}

\noindent It follows a minimal description of $P(K_n)$ ($n\geq 3$).

\begin{cor}
A minimal description of $P(K_n)$ ($n\geq 3$) is the set of $x\in \mathbb{R}^E$ satisfying (1), (3) and constraints (4) with $2\leq |U|\leq \lfloor{\frac{n}{2}}\rfloor$.
\end{cor}

\noindent A natural question is: what about any graph? In this paper we generalize this result by giving a minimal description of $P(G)$ in the general case and by presenting a short proof of it. Another motivation is that "Facets are strongest cutting planes in an integer programming sense (see \cite{Garfinkel and Nemhauser 1972}) and it is thus natural to expect that such inequalities are of substantial computational value in the numerical solution of this hard combinatorial optimization problem" \cite{Grotschel and Padberg 1979a}.
\newline Since one of the motivations for this question is to optimize over $P(G)$, then we can suppose that $G$ is a 2-connected graph without loops (deletion of the loops), bridges ($P(G)=\O$), parallel (deletion of one parallel edge) or series edges (contraction of one series edge) because any of these cases can be reduced to our case by using the indicated operation. We use the notations: $n=|V|$, $m=|E|$, and $n_H=|V(H)|$, $m_H=|E(H)|$ for any subgraph or any subset of edges $H$ of $G$. Moreover for $y\in \mathbb{R}^E$ and $i\in \mathbb{Z}$, $support_i(y)=\{ e\in E$ such that $y(e)=i\}$. For $k\in \mathbb{Z}_+$ and $\mathcal X=\{ x_1, ..., x_k\}\subset \mathbb{R}^E$, $\mathcal X$ satisfy the intersection condition (IC) (respectively union condition (UC)) if $\bigcap_{i=1}^{k} support_1 (x_i)=\O$ (respectively $\bigcup_{i=1}^{k} support_0 (x_i)=E$), which is equivalent to the following: for any $e\in E$, there exists $i\in \{ 1, ..., k\}$ such that $x_i(e)\neq 1$ (respectively $x_i(e)=0$). It follows that if $\mathcal X\supseteq \mathcal X'$ (respectively $\mathcal X\subseteq \mathcal X'$) and $\mathcal X$ satisfy IC (respectively UC) then $\mathcal X'$ too. IC and UC are equivalent when the $x_i$'s are $\{ 0, 1\}$-vectors.
\newline Let $M$ be a matroid defined on a finite set $E$. $\mathcal {B}(M)$, $M^*$, the functions $r$ and $r^*$, and $K(M)$ are, respectively, the class of bases, the dual matroid, the rank and the dual rank  functions, and the bases polytope of $M$.
Suppose that $M$ (and $M^*$) is 2-connected. A subset $L\subset E$ is called a locked subset of $M$ if $M|L$ and $M^*|(E\backslash L)$ are 2-connected, and their corresponding ranks are at least 2, i.e., $min\{r(L), r^*(E\backslash L)\} \geq 2$. It is not difficult to see that if $L$ is locked then both $L$ and $E\backslash L$ are closed, respectively, in $M$ and $M^*$ (That is why we call it locked). Locked subsets were introduced to solve many combinatorial problems in matroids \cite{Chaourar 2002, Chaourar 2008, Chaourar 2011, Chaourar 2017, Chaourar 2018}.
\newline By analogy, $K(G)=K(M(G))$ and a locked subgraph $H$ of $G$ is a subgraph for which $E(H)$ is a locked subset of $M(G)$.
\newline We can write a minimal description of $K(G)$ as a corollary of a similar result for $K(M)$ \cite{Chaourar 2018, Feichtner and Sturmfels 2005, Fujishige 1984}.

\begin{thm} A minimal description of $K(G)$ is the set of all $x\in \mathbb{R}^E$ satisfying the following inequalities:
$$  x(e)\geq 0  \> for\>any\> edge\> e\in E	                                     \eqno   (5)$$
$$	x(e) \leq 1 \> for\>any\> edge\> e\in E			                             \eqno 	 (6)$$
$$	x(E(H)) \leq n_H-1	 \> for\> any\> locked\> subgraph\> H\>of\> G	         \eqno	 (7)$$
$$  x(E)=n-1                                                                     \eqno   (8)$$
\end{thm}

\noindent For an induced subgraph $H$ of $G$, $\overline{H}=(V(E\backslash E(H)), E\backslash E(H))$ is called the complementary subgraph of $H$ in $G$. Moreover, for any $F\subseteq E$, $H\cup F$ is the subgraph $(V(H)\cup V(F), E(H)\cup F)$.
\newline A collection $\mathcal C$ of subsets of a nonempty set $X$ is called laminar if for all $T,U\in \mathcal C$: $T\subseteq U$ or $U\subseteq T$ or $T\cap U=\O$. There is the following upper bound on the size of a laminar family \cite{Schrijver 2004}:
\begin{thm}
If $\mathcal C$ is laminar and $X\notin \mathcal C$, then $|\mathcal C|\leq 2|X|-1$.
\end{thm}
\noindent The remainder of the paper is organized as follows: in section 2, we give a minimal description of $P(G)$ in the general case, then we conclude in section 3.
\section{Facets of the subtours elimination polytope}

First, we need the following lemma.

\begin{lem} Let $H$ be a 2-connected subgraph of $G$, and $\{ L_1, L_2\}$ be a partition of $E(\overline{H})$. Then $\overline{H}$ is connected if and only if $n_H+n<n_{H\cup L_1}+n_{H\cup L_2}$.
\end{lem}
\bpf $n_H+n<n_{H\cup L_1}+n_{H\cup L_2}$ is equivalent to: $n_H+n_H+|V(\overline{H})|-|V(H)\cap V(\overline{H})|<n_H+n_{L_1}-|V(H)\cap V(L_1)|+n_H+n_{L_2}-|V(H)\cap V(L_2)|$, i.e., $2n_H+|V(\overline{H})|-|V(H)\cap V(\overline{H})|<2n_H+n_{L_1}+n_{L_2}-|V(H)\cap V(L_1)|-|V(H)\cap V(L_2)|=2n_H+n_{L_1}+n_{L_2}-|V(H)\cap V(\overline{H})|-|V(H)\cap V(L_1)\cap V(L_2)|$, i.e., $n_{L_1}+n_{L_2}-|V(L_1)\cap V(L_2)|=|V(\overline{H})|< n_{L_1}+n_{L_2}-|V(H)\cap V(L_1)\cap V(L_2)|$, i.e., $|V(L_1)\cap V(L_2)|>|V(H)\cap V(L_1)\cap V(L_2)|\geq 0$, or  $|V(L_1)\cap V(L_2)|\geq 1$, which means that $\overline{H}$ is connected.
\epf

\noindent Now we can characterize locked subgraphs by means of graphs terminology.

\begin{lem} $H$ is a locked subgraph of $G$ if and only if $H$ is an induced and 2-connected subgraph such that $3\leq n_H\leq n-1$, and $\overline{H}$ is a connected subgraph.
\end{lem}
\bpf
It is not difficult to see that M(H) is closed and 2-connected in M(G) if and only if $H$ is an induced and 2-connected subgraph of $G$.
\noindent Now, suppose that $E(H)$ is closed and 2-connected, and $E(\overline{H})$ is 2-connected in the dual matroid $M^*(G)$. Let $\{ L_1, L_2\}$ be a partition of $E(\overline{H})$. It follows that $r^*(E(\overline{H}))<r^*(L_1)+r^*(L_2)$, i.e., $|E(\overline{H}|-r(E)+r(E(H))<|L_1|+|L_2|-2r(E)+r(E(H\cup L_1)|+r(E(H\cup L_2)|$. In other words, $r(E(H))+r(E)=r(E(H)\cup L_1)+r(E(H)\cup L_2)$, which is equivalent to: $n_H-1+n-1<n_{H\cup L_1}-1+n_{H\cup L_2}-2$, i.e., $\overline{H}$ is connected according to the previous lemma.
\newline Let check the condition: $min\{r(E(H)), r^*(E(\overline{H}))\} \geq 2$. Since $r(E(H))=n_H-1$ then we have $r(E(H))\geq 2$ if and only if $n_H\geq 3$. Moreover, $r^*(E(\overline{H}))=m_{\overline{H}}+r(E(H))-r(E)=m_{\overline{H}}+n_H-n$ then we have $r^*(E(\overline{H}))\geq 2$ if and only if $n_H\geq 2+n-m_{\overline{H}}$, i.e., $|V(H)\cap V(\overline{H})|\geq 2+n_{\overline{H}}-m_{\overline{H}}$. But, since $G$ is 2-connected and $\overline{H}$ is connected, either $|V(H)\cap V(\overline{H})|\geq 3$, or $|V(H)\cap V(\overline{H})|=2$ and $n_{\overline{H}}\leq m_{\overline{H}}$. In both cases, the later inequality is verified and we do not need to mention it.
\newline Furthermore, $M(H)$ is closed and distinct from $E$, i.e., $r(E(H))\leq r(E)-1$, is equivalent to: $n_H\leq n-1$.
\epf

\noindent We have then the following refined description of the subtours elimination polytope.

\begin{lem} Let $v_0\in V$. $P(G)$ is the set of all $x\in \mathbb{R}^E$ satisfying (1), (2), and the following constraints:
$$	x(\delta (v))=2 \> for\>any\> vertex\> v\in V\backslash \{ v_0\}               \eqno	 (8)  $$
$$	x(\delta (U))\geq 2	 \> for\> any\> locked\> subgraph\> G(U)\> of\> G 	       \eqno	 (9)  $$
$$  x(E)=n                                                                         \eqno    (10) $$						
\end{lem}
\bpf
It is not difficult to see that (3) is equivalent to (8) and (10).
\newline We will prove now the equivalence of a constraint (9). It is clear that $G(U)$ is connected because of (9). Suppose that $G(U)$ is not 2-connected, then there exist two subsets $U_1$ and $U_2$ of $U$ such that $G(U)=G(U_1)\oplus G(U_2)$. Let $u$ be the unique vertex in $U_1\cap U_2$. It follows that $x(\delta (U))=x(\delta (U_1))+x(\delta (U_2))-x(\delta (u))\geq 2+2-2=2$ and the constraint (4) is redundant for $U$. Suppose now that $\overline{G(U)}$ is not connected. It follows that there exists a partition $\{ U_1, U_2\}$ of $U$ such that $x(\delta (U))=x(\delta (U_1))+x(\delta (U_2))$. It follows that (4) is redundant for $U$. The constraints (8) and (10) imply the redundancy of a constraint (4) for $|U|=1$. Suppose now that $|U|=2$, $U=\{ u, w\}$, and let $e=uw\in E$. If follows that $x(\delta (U))=x(\delta (u))+x(\delta (v))-2x(e)\geq 2+2-2=2$, so constraint (4) is redundant for $U$. At last, since $U$ is a proper subset of $V$ then $|U|\leq n-1$.
\epf

\noindent Let $Q(G)$ be the set of $x\in \mathbb{R}^E$ satisfying (1), (2), (10) and:
$$	x(H)\leq n_H-1	 \> for\> any\>  locked\> subgraph\> H\> of\> G 	    \eqno     (11)$$
\noindent Let $Q'(G)$ be the set of $x\in \mathbb{R}^E$ satisfying (1), (2), (9), (10) and:
$$	x(\delta (v))\geq 2 \> for\>any\> vertex\> v\in V                       \eqno	  (12)  $$
\newline It follows that:

\begin{lem}
(1) $P(G)\subseteq Q(G)\subseteq Q'(G)$. (2) Moreover, the extreme integer points of $Q(G)$ are exactly the Hamilton circuits of $G$.
\end{lem}
\bpf
(1) Since $2x(E(U))+x(\delta (U))=\sum_{u\in U} x(\delta (u)) =2|U|$, then a constraint (11) with $|V(H)|=|V|-1$ implies a constraint (12). It is not difficult to see, similarly as in the previous lemma proof, that constraint (11) implies the following constraint:
$$	x(E(U))\leq |E(U)|-1	 \> for\> any\>  proper\> subset\> U\> of\> V 	    \eqno     (11)$$
Thus, constraint (11) implies constraint (9), for a similar reason as for constraint (12), and $Q(G)\subseteq Q'(G)$. According to the previous lemma, $P(G)\subseteq Q(G)$.
\newline (2) It is clear that Hamilton circuits of $G$ are integer extreme points of $Q(G)$ because they satisfy all constraints of $Q(G)$ and $m-1$ independent constraints of type (1) and (2). Now, let $x\in Q(G)$ be an integer extreme point. It follows that $x$ is a $\{ 0, 1\}$-vector in $Q'(G)$ and $x(E)=n$, thus $x$ is a Hamilton circuit.
\epf

So optimizing over $Q(G)$ will give a lower (respectively an upper) bound for TSP (respectively MaxTSP).
\newline We need the following lemma for the coming theorem.

\begin{lem}
Let $x\in Q(G)$ be an extreme point and $\mathcal L_x$ be an independent system of tight constraints of type (11) verified by $x$. Then we can choose a laminar system $L_x=\{U\subset V$ such that $x(E(U))=|U|-1\}$ equivalent to $\mathcal L_x$ and $|L_x|=|\mathcal L_x|$.
\end{lem}
\bpf
Since $y(E(U))+y(E(W))=y(E(U\cap W))+y(E(U\cup W))$ for any $U, W\subset V$ and $y\in \mathbb{R}_+^E$, then $U$ and $W$ define tight constraints of type (11) if and only if $U\cap W$ and $U\cup W$ too. Applying intersections and unions to subsets of $\mathcal L_x$ will transform it into $L_x$.
\epf

\noindent We have the following characterization of extreme points of $Q(G)$.

\begin{thm}
(1) For any extreme point $x\in Q(G)$, there are $x_i\in K(G)$, $i=1, ..., n$ such that $\{ x_1, ..., x_n\}$ satisfy IC and $x=\frac{1}{n-1} \sum_{i=1}^{n} x_i$.
\newline (2) For any $x_i\in K(G)$, $i=1, ..., n$, if $\mathcal X=\{ x_1, ..., x_n\}$ satisfy UC then $\frac{1}{n-1} \sum_{i=1}^{n} x_i\in Q(G)$.
\end{thm}
\bpf
(1) {\bf Case 1:} $x$ is integer. It follows that $x$ is a Hamilton circuit and the $x_i$'s are the $n$ Hamilton paths contained in $x$.
\newline {\bf Case 2:} $x$ is not integer. Let $A=\{ e\in E$ such that $x(e)\geq \frac{1}{n-1}\}$, $E_1 =support_1(x)$, $E_x =E\backslash support_0(x)$, $G_x=(V, A)$, and $c$ is the number of connected components of $G_x$. We can suppose that $|E_1|\leq n-2$ because otherwise any $T\subseteq E_1$ with $|T|=n-1$ is a spanning tree and $y=x-\frac{1}{n-1} T\in K(G)$.
\newline {\bf Claim 1:} $A\neq \O$.
\newline By contradiction, suppose that $x(e)< \frac{1}{n-1}$ for any edge $e\in E$. Thus, $n=x(E)<\frac{m}{n-1}\leq \frac{n(n-1)}{2} \frac{1}{n-1}=\frac{n}{2}$, a contradiction.
\newline {\bf Claim 2:} $|A|\geq n$.
\newline By contradiction, suppose that $|A|\leq n-1$. Thus, there exists a vertex $v$ such that $|A\cap \delta (v)|\leq 1$ because, otherwise, $|A|=\frac{1}{2} \sum_{v\in V} |A\cap \delta (v)|\geq \frac{1}{2} 2n=n$, and we are done. Since $x(\delta (v)\backslash A)+x(A\cap \delta (v))=x(\delta (v))\geq 2$, because $Q(G)\subseteq Q'(G)$, and $x(A\cap \delta (v))\leq |A\cap \delta (v)|\leq 1$, because of constraint (2), then $\frac{|\delta (v)|}{n-1}>x(\delta (v)\backslash A)\geq 1$, i.e., $|\delta (v)|> n-1$, a contradiction.
\newline {\bf Claim 3:} $G_x$ is a connected subgraph.
\newline By contradiction, suppose that $c\geq 2$ and let $(V_i, E_i)$, $i=1, ..., c$, be the connected components of $G_x$. For any $i=1, ..., c$, $2\leq x(\delta (V_i))<\frac{|\delta (V_i)|}{n-1}$. It follows that $|\delta (V_i)|\geq 2n-1$. Thus, $|E_x\backslash A|\geq \frac{1}{2} \sum_{i=1}^{c} |\delta (V_i)|\geq cn-\frac{c}{2}$, and then, according to Claim 2, $|E_x|=|A|+|E_x\backslash A|\geq (c+1)n-\frac{c}{2}$. Since $x$ is an extreme point of $Q(G)$ then it should verify at least $\ell$ tight constraints of type (11) that we can choose laminar according to the previous lemma. Thus $\ell \geq |E_x|-|E_1|-1\geq (c+1)n-\frac{c}{2}-n+2-1=cn-\frac{c}{2}+1=\frac{c(2n-1)+2}{2}\geq \frac{2(2n-1)+2}{2}=2n$. But according to Theorem 1.4, $\ell \leq 2n-1$, a contradiction.
\newline {\bf Claim 4:} $E_1$ does not contain a circuit.
\newline $|E_1|=x(E_1)\leq |V(E_1)|-1<|V(E_1)|$. Thus, $E_1$ is not a circuit. Actually the later inequality is true for any subset of $E_1$ and we are done.
\newline It follows that there exists a spanning tree $T\supset E_1$ in $G_x$ such that $y=x-\frac{1}{n-1} T\in K(G)$ because of the followings. $y(e)\geq 0$ because $T\subseteq A$. $y(e)\leq x(e)\leq 1$, for any edge $e\in E$. Similarly, $y(E(H))\leq x(E(H))\leq n_H-1$, for any locked subgraph $H$ of $G$. At last, $y(E)=x(E)-\frac{1}{n-1} T(E)=n-\frac{n-1}{n-1}=n-1$.
\newline In this case, $y$ can be written as a convex combination of spanning trees as follows: $y=\sum_{j=1}^{t} \lambda_j T_j=\frac{1}{n-1}\sum_{j=1}^{t} (n-1)\lambda_j T_j=\frac{1}{n-1}\sum_{j=1}^{t'} \lambda'_j T'_j$ with $\sum_{j=1}^{t'} \lambda'_j =\sum_{j=1}^{t} (n-1)\lambda_j =(n-1)sum_{j=1}^{t} \lambda_j =n-1$ and $0\leq \lambda'_j\leq 1$. So this sum can be partitioned into $n-1$ convex combinations of spanning trees: $y=\frac{1}{n-1}\sum_{i=1}^{n-1} \sum_{j=1}^{t_i} \mu_{ij} T_{ij}=\frac{1}{n-1}\sum_{i=1}^{n-1} x_i$ with $x_i\in K(G)$, $i=1, ..., n-1$.
\newline Now, since $T\supset E_1$, then $y(e)<1$ for any edge $e\in E$, i.e., $\frac{1}{n-1}\sum_{i=1}^{n-1} x_i(e)<1$. Let $i_0\in \{ 1, ..., n-1\}$ such that $x_{i_0}(e)=Min_{1\leq i\leq n-1} \{ x_i(e)\}$. Thus, $x_{i_0}(e)=\frac{1}{n-1}\sum_{i=1}^{n-1} x_{i_0}(e)\leq \frac{1}{n-1}\sum_{i=1}^{n-1} x_i(e)<1$,, i.e., $x_{i_0}(e)\neq 1$ and we are done.
\newline (2) First we prove (2) when the $x_i$'s are spanning trees. Let $x=\frac{1}{n-1} \sum_{i=1}^{n} T_i$. It is clear that $x$ satisfy (1) and (10). UC for $\mathcal T$ is equivalent to the following. For any edge $e\in E$, there exists $i\in \{1, ..., n\}$ such that $T_i(e)=0$. Thus, $\sum_{i=1}^{n} T_i(e)\leq n-1$, and then $x(e)\leq \frac{n-1}{n-1}=1$. Now let $H$ be a locked subgraph of $G$ and $E_i(H)=\{e\in E(H)$ such that $T_i(e)=0\}$. It follows that $T_i(E(H))=n_H-1-|E_i(H)|$, i.e., $\sum_{i=1}^{n} T_i(E(H))=n(n_H-1)-\sum_{i=1}^{n} |E_i(H)|\leq n(n_H-1)-m_H$, because $E(H)=\bigcup_{i=1}^{n} E_i(H)$ (UC), i.e., $m_H\leq \sum_{i=1}^{n} |E_i(H)|$. Thus, $\sum_{i=1}^{n} T_i(E(H))\leq n(n_H-1)-(n_H-1)=(n-1)(n_H-1)$, and we are done.
\newline Now let $x=\frac{1}{n-1} \sum_{i=1}^{n} x_i$ with $\mathcal X=\{ x_1, ..., x_n\}\subset K(G)$ satisfy UC. Since $x_i\in K(G)$ then $x_i$ can be written as a convex combination of spanning trees: $x_i=\sum_{j=1}^{k_i} \lambda_{i, j} T_{i, j}$. UC for $\mathcal X$ is equivalent to the following. For any edge $e\in E$, there exists $i\in \{1, ..., n\}$ such that $x_i(e)=0$, i.e., $T_{ij}(e)=0$ for all $j=1, ..., k_i$. Thus, $x=\frac{1}{n-1} \sum_{i=1}^{n} x_i=\frac{1}{n-1} \sum_{i=1}^{n} \sum_{j=1}^{k_i} \lambda_{i, j} T_{i, j}=\sum_{j=1}^{k} \lambda'_j \frac{1}{n-1} \sum_{i=1}^{n} T_{i, j}=\sum_{j=1}^{k} \lambda'_{i, j} y_j$ (e.g. $\lambda'_1=Min_{i=1}^{n} \{ \lambda_{i, 1}$ and so on). But $y_j\in Q(G)$ for all $j=1, ..., m$ (because for any edge $e\in E$, there exists a spanning tree $T_{i, j}$ in the expression of $y_j$ such that $T_{ij}(e)=0$, i.e., $y_j(e)\neq 1$) and $Q(G)$ is convex, then $x\in Q(G)$.
\epf

For the coming results, we need the following notations.
\begin{itemize}
\item For a polytope $A$, $pA$ is the set of extreme points of $A$.
\item Let $k\in \mathbb{Z}_+$. For $A\subseteq \mathbb{R}^k$, $PA$ denotes the smallest convex set of $\mathbb{R}^k$ containing $A$.
\item For $A, B\subseteq K(G)$ and $\lambda\in \mathbb{R}$, $A+B=\{ x+y$ such that $(x, y)\in A\times B\}$, and $\lambda A=\{ \lambda x$ such that $x\in A\}$,
\item $\sum_{i=1}^{k} A_i=A_1+ ...+A_k$.
\item $\prescript{k}{}{A}=\sum_{i=1}^{k} A_i$ with $A_i=A$, for $i=1, ..., k$.
\item Let $H$ be a locked subgraph of $G$. $F_H=\{ x\in Q(G)$ such that $x(E(H))=n_H-1\}$,
\item $F_i=\{x\in K(G)$ such that $x(E(H))=n_H-1-i\}$,
\item $K_i=\{ x\in K(G)$ such that $x(E(H))\leq n_H-1-i\}$, $i\in \mathbb{Q}$,
\item $\Omega =\{(k_1, ..., k_{n_H-1})\in \mathbb{Z}_+^{n_H-1}$ such that $\sum_{i=0}^{n_H-2} k_i (n_H-1-i)=(n_H-1)(n-1)$ and $\sum_{i=0}^{n_H-2} k_i\leq n\}$,
\item $F'_{(k_0, ..., k_{n_H-2})}=\sum_{i=0}^{n_H-2}\prescript{k_i}{}F_i$,
\item $F'_\lambda=\bigcup_{(k_0, ..., k_{n_H-2})\in \Omega} \lambda F'_{(k_0, ..., k_{n_H-2})}$,
\item $K'=\prescript{n}{}(K(G))$.
\item $\mathcal T(G)$ be the class of spanning trees of $G$,
\item $\mathcal T_i(G)=\{ T\in \mathcal T(G)$ such that $|T\cap E(H)|=n_H-1-i\}$,
\item $Q''=\bigcup_{i=0}^{n_H-1} (\frac{1}{n-1} \mathcal T_i+K_{\frac{n_H-1-i}{n-1}})$,
\item $K''=\frac{1}{n-1} \mathcal T+K(G)$
\item For $X\subseteq \prescript{k}{} (\mathbb{R}^E)$, $X(OC)=\{ x\in X$ such that $x$ satisfy OC$\}$ with OC$\in\{$ IC, UC$\}$.
\end{itemize}
Note that we have removed from some notations the index $H$ for more readability.
\newline A second corollary of Theorem 2.6 follows.

\begin{cor}
$PK'_{\frac{1}{n-1}}(UC)\subseteq Q(G)\subseteq PK'_{\frac{1}{n-1}}(IC)$.
\end{cor}
\bpf
$K'_{\frac{1}{n-1}}(UC)\subseteq Q(G)$ and $pQ(G)\subseteq K'_{\frac{1}{n-1}}(IC)$ according to the Theorem 2.6. By applying $P$ to both inclusions, we get the result.
\epf

We need the following lemma.

\begin{lem}
If $(k_1, ..., k_{n_H-2})\in \Omega$ then $k_0\geq n-n_H+1(\geq 2)$. 
\end{lem}
\bpf
We have $\sum_{i=0}^{n_H-2} k_i (n_H-1-i)=(n_H-1)(n-1)$, i.e., $\sum_{i=1}^{n_H-2} k_i (n_H-1-i)=(n_H-1)(n-1-k_0)$. It follows that: $(n_H-2)(n-k_0)\geq (n_H-2)\sum_{i=1}^{n_H-2} k_i\geq \sum_{i=1}^{n_H-2} (n_H-1-i)k_i=(n_H-1)(n-1-k_0)$, i.e., $(n_H-1)(n-k_0)-(n-k_0)\geq (n_H-1)(n-k_0)-(n_H-1)$. Thus, $k_0\geq n-n_H+1$.
\epf

We use the following properties of the dimension.

\begin{lem} \
\begin{enumerate}
\item $Max\{ dim(A), dim(B)\}\leq dim(A+B)\leq dim(A)+dim(B)$.
\item $dim(\lambda A)=dim(A)$ if $\lambda \neq 0$.
\item $dim(A)=0$ if $A$ is a finite set.
\item $dim(F_1)=m-2$ and $dim(K(G))=m-1$.
\end{enumerate}
\end{lem}

\noindent Now, we can state our main result.

\begin{thm}
The constraints (1)-(2) and (10)-(11) give a minimal description of $Q(G)$.
\end{thm}
\bpf
It is not difficult to prove that constraints (1)-(2) and (10) are irredundant. We need only to prove that $dim(F_H)=m-2$.
\newline We have $PF'_{\frac{1}{n-1}}(UC)\subseteq F_H$ because $F'_{\frac{1}{n-1}}(UC)\subseteq F_H$.
\newline It follows that $dim(F'_{\frac{1}{n-1}}(UC))\leq dim(F_H)$. Moreover, $m-2=dim(F_1)\leq dim(F'_{\frac{1}{n-1}}(UC))$ (because of Lemma 2.8) and $dim(F_H)\leq dim(Q(G)=m-1$. Thus, $m-2\leq dim(F'_{\frac{1}{n-1}}(UC))\leq dim(F_H)\leq m-1$.
\newline By using the proof of Theorem 2.6, we can see that: $F_H\subseteq \bigcup_{i=0}^{n_H-1} (\frac{1}{n-1} \mathcal T_i+F_{\frac{n_H-1-i}{n-1}})$. It follows that $F_H\subseteq Q''$, and, since $dim(F_{\frac{n_H-1-i}{n-1}})<dim(K_{\frac{n_H-1-i}{n-1}})$, then $dim(F_H)<dim(Q'')$. In the other hand, $Q''\subseteq K''$, i.e., $dim(Q'')\leq dim(K'')$. At last $dim(K'')\leq dim(\mathcal T(G))+dim(K(G)=m-1$ and we are done.
\epf

\noindent And we have the following consequence for $P(G)$.

\begin{cor}
The constraints (1), (3) and (9) with $3\leq |U|\leq n-2$ give a minimal description of $P(G)$.
\end{cor}

\section{Conclusion}

We have characterized all facets of $P(G)$ and given a minimal description of it. Future investigations can be improving the fractional ratio of integrity of $P(G)$ and providing a combinatorial algorithm for optimizing on $P(G)$.

\end{document}